\documentclass[11pt]{amsart}

\newcommand{\E}{\mathbb{E}}
\newcommand{\Z}{\mathbb{Z}}
\newcommand{\N}{\mathbb{N}}
\newtheorem{theorem}{Theorem}

\newtheorem*{proposition*}{Proposition}

\newtheorem{lemma}{Lemma}

\theoremstyle{definition}
\newtheorem{definition}{Definition}

\newcommand{\norm}[1]{\left\Vert #1\right\Vert}

\begin{document}
\title[Multiple recurrence and convergence]{Multiple recurrence and convergence for sequences related to the prime numbers}
\author{Nikos Frantzikinakis, Bernard Host, and Bryna Kra}

\begin{abstract}
For any measure preserving system $(X,\mathcal{X},\mu,T)$ and
$A\in\mathcal{X}$ with $\mu(A)>0$, we show that there exist
infinitely many primes $p$ such that $\mu\bigl(A\cap
T^{-(p-1)}A\cap T^{-2(p-1)}A\bigr) > 0$ (the same holds with $p-1$
replaced by $p+1$). Furthermore, we show the existence of the
limit in $L^2(\mu)$ of the associated ergodic average over the
primes. A key ingredient is a recent result of Green and Tao on
the von Mangoldt function. A combinatorial consequence is that
every subset of the integers with positive upper density contains
an arithmetic progression of length three and common difference of
the form $p-1$ (or $p+1$) for some prime $p$.
\end{abstract}

\address{Institute for Advanced Study, 1 Einstein drive,
Princeton, NJ 08540}
\address{\'Equipe d'analyse et de math\'ematiques appliqu\'ees,
Universit\'e de Marne la Vall\'ee, 77454 Marne la Vall\'ee Cedex,
France}
\address{Department of Mathematics, Northwestern University,
2033 Sheridan Road,  Evanston, IL 60208-2730, USA}
\email{nikos@ias.edu} \email{Bernard.Host@univ-mlv.fr}
\email{kra@math.northwestern.edu}
\thanks{The first author acknowledges the support of
NSF grant DMS-0111298 and the third author of NSF grant
DMS-0555250.}



 \maketitle

\section{Introduction}
\subsection{Results about the primes}
The \emph{von Mangoldt function} $\Lambda\colon
\mathbb{Z}\to\mathbb{R}$, defined by
$$
\Lambda(n) =
\begin{cases}
\log p & \text{ if } n = p^m \text{ for some } m\in\mathbb{N}
\text{ and } p\in\mathbb{P}\\
0  & \text{ otherwise}
\end{cases}
$$
plays a major role in understanding the
distribution of the prime numbers.
The classical circle method is the study of the Fourier transform of
$\Lambda$ restricted to an interval $[1,N]$,
that is, evaluating
the correlation of $\Lambda$ with complex exponentials on this
interval.


Green and Tao (\cite{GT}, \cite{GT2}, \cite{GT3}) generalize this
method by estimating the correlations of $\Lambda$ with
``nilsequences'' and then deduce an estimate for the third Gowers
norm of a modification of $\Lambda$ (a precise definition of this
norm is given in Section~\ref{sec:nt}). This result is the base of
our study; we need some notation to state it explicitly.

Let $\mathbb{P}$ be the set of prime numbers.
If $w$ is a positive integer and $r\in\mathbb{Z}$, setting
$$
W=\prod_{p\in \mathbb{P}, p < w}p\;,
$$
for $n\in \N$ we define
$$
\tilde{\Lambda}_{w,r}(n)=\frac{\phi(W)}{W}\cdot \Lambda(Wn+r)\ ,
$$
where $\phi$ is the Euler function.
It is easier to work with this modified von Mangoldt function, as it is
better distributed over congruence classes than the original
function.

If $N$ is a prime number we identify $[0,N-1]$ with $\Z/N\Z$ in
the natural way and consider the restriction
$\tilde{\Lambda}_{w,N,r}(n)$ of $\tilde{\Lambda}_{w,r}(n)$ to this
interval as a function on $\Z/N\Z$.  It follows immediately from
the results in \cite{GT3} that (notation explained in
Section~\ref{sec:nt}):

\begin{theorem}[\bf Green and Tao~\cite{GT3}]\label{C:GT}
For $r\in\mathbb{Z}$ with $(r,W)=1$ we have
$$
\norm{ \bigl(\tilde{\Lambda}_{w,N, r }(n)-1\bigr)\cdot {\bf
1}_{[0,[N/3])}(n)}_{U_3(\Z/N\Z)}=o_{N\to\infty;
w}(1)+o_{w\to\infty}(1).
$$
\end{theorem}

We use this result to derive several results in ergodic theory, on
recurrence and convergence properties of some sequences related to
the primes, and use them to deduce a combinatorial consequence.

\subsection{Multiple recurrence along primes}
The connection between additive combinatorics and ergodic theory
began with Szemer\'edi's celebrated theorem~\cite{S}, which states
that a subset of the  integers with positive upper density
contains arbitrarily long arithmetic progressions.
Furstenberg~\cite{F1} found an equivalent formulation of
Szemer\'edi's Theorem in terms of multiple recurrence.  He proved
this ``correspondence principle'' and showed that if
$(X,\mathcal{X},\mu,T)$ is a measure preserving system\footnote{A
measure preserving system is a quadruple  $(X,\mathcal{X},\mu,
T)$, where $(X,\mathcal{X},\mu)$  is a probability space and
$T\colon X\to X$ is a measurable map such that
$\mu(T^{-1}A)=\mu(A)$ for all $A\in\mathcal{X}$. Without loss of
generality we can assume that the probability space is Lebesgue.},
$A\in\mathcal{X}$ has positive measure, and $k$ is a positive
integer, then there are infinitely many positive integers $n$ such
that
\begin{equation}
\label{eq:multrec}
\mu\bigl(A\cap T^{-n}A\cap
T^{-2n}A\cap\ldots\cap T^{-kn}A\bigr) > 0 \ .
\end{equation}

A natural question is what restrictions can be placed on $n$ such
that  the measure of the intersection in~\eqref{eq:multrec}
remains positive.  This motivates the following definition:

\begin{definition}
Given an integer $k\geq 1$, $S\subset \mathbb{N}$ is a \emph{set
of $k$-recurrence} if for every measure preserving system
$(X,\mathcal{X},\mu,T)$ and $A\in\mathcal{X}$ with $\mu(A)>0$,
there exist infinitely many $n\in S$ such that
$$
  \mu\bigl(A\cap T^{-n}A \cap T^{-2n}\cap \cdots \cap T^{-kn}A\bigr)>0\ .
$$
\end{definition}

Via Furstenberg's correspondence principle, the equivalent
combinatorial formulation of this definition is classical:
\begin{proposition*}
Let $k\geq 1$ be an integer.  The set $S\subset \mathbb{N}$ is  a
set of $k$-recurrence if and only if every subset $A$ of integers
with positive upper density contains at least one arithmetic
progression of length $k+1$ and common difference in $S$.
\end{proposition*}
A set of $k$-recurrence is also known in the number theory
literature as a \emph{$k$-intersective set}.

There are many known examples of sets of $1$-recurrence.  For
example, one can take all multiples of a fixed number
or $S-S$ for any infinite set $S\subset\mathbb{N}$.
Furstenberg~\cite{F2} and
S\'ark\"ozy~\cite{Sa1} independently
showed that values of a polynomial, meaning to $\{q(n)\colon n\in
\mathbb{N}\}$ where $q(n)$ is an integer valued polynomial with
$q(0)=0$, form a set of $1$-recurrence.
Another interesting example
comes from the primes.
S\'ark\"ozy~\cite{Sa2} showed that the set of
shifted primes $\mathbb{P}-1$ (and the set $\mathbb{P}+1$)
form a set of $1$-recurrence.
Wierdl~\cite{W} reproved this result using methods
from ergodic theory.

For higher $k$, restricting the exponents in~\eqref{eq:multrec} is
more difficult. It is an immediate consequence of Szemer\'edi's
theorem that one can take the set of all multiples of a fixed
nonzero integer. Bergelson and Leibman~\cite{BL} showed that
polynomial values $\{q(n)\colon n\in \mathbb{N}\}$ where $q(n)$ is
an integer valued polynomial with $q(0)= 0$ are sets of
$k$-recurrence for all $k\geq 1$. On the other hand,
Furstenberg~\cite{F2} constructed an infinite set $S$ whose set of
differences $S-S$  is not a set of $2$-recurrence. For general
$k$, a set of $k$-recurrence but not $(k+1)$-recurrence was
constructed in~\cite{FLW}.

The recent results of Green and Tao on the von Mangoldt
function allow us to show $2$-recurrence for the primes:


\begin{theorem}\label{T:recurrence1}
\label{th:recurrence}
Let $(X, \mathcal{X}, \mu, T)$ be a measure preserving system and
let $A\in\mathcal{X}$ have positive measure.
There exist infinitely many $p\in\mathbb{P}$ such
that
\begin{equation*}
\mu\bigl(A\cap T^{-(p-1)}A\cap
T^{-2(p-1)}A\bigr) > 0
\end{equation*}
and there exist infinitely many $p\in\mathbb{P}$ such
that
\begin{equation*}
\mu\bigl(A\cap T^{-(p+1)}A\cap
T^{-2(p+1)}A\bigr) > 0 \ .
\end{equation*}
\end{theorem}

It is easy to check that the set $\mathbb{P}+r$ is not even a set
of $1$-recurrence for $r\in \mathbb{Z}\setminus \{-1,1\}$, by
considering a periodic system with period $r$ if $r\neq 0$ and
with period $2$ if $r=0$. A related topological version of this
question was posed by Brown, Graham, and Landman~\cite{BGL}.

A corollary of Furstenberg's correspondence principle is
that any set of integers with positive upper density
contains arithmetic progressions of length $3$ with common
difference in $\mathbb{P}-1$ (and also in $\mathbb{P}+1$).

The higher order statement of Theorem~\ref{C:GT} was conjectured
to be true in \cite{GT3}, if it holds then the obvious
generalizations of our proofs show that $\mathbb{P}-1$ and
$\mathbb{P}+1$ are sets $k$-recurrence for every $k\in\N$. This
and other generalizations are discussed in
Section~\ref{sec:further}.

\subsection{Convergence of averages along primes}
A closely related question is the convergence of the associated
multiple ergodic averages.  In his proof of Szemer\'edi's Theorem,
Furstenberg proved more than the measure of the intersection in
\eqref{eq:multrec} is positive.
He showed that given a  measure
preserving system $(X, {\mathcal X}, \mu, T)$ and $A\in{\mathcal
X}$ of positive measure, we have
\begin{equation}
\label{eq:AP} \liminf_{N\to\infty} \, \frac{1}{N}\sum_{0\leq n<N}
\mu\bigl(A\cap T^{-n}A\cap T^{-2n}A\cap\ldots\cap T^{-kn}A\bigr)
> 0 \ .
\end{equation}

A natural question is the existence of the limit of the associated
ergodic average. More generally, given a measure preserving system
$(X,\mathcal{X},\mu,T)$, functions $f_1,$ $f_2,$ $\ldots,$ $f_k\in
L^\infty(\mu)$, and an increasing sequence of integers
$\{s_n\}_{n\in \N}$, we can study the existence of the limit of
$$
\frac{1}{N}\sum_{0\leq n<N}(T^{s_n}f_1\cdot
T^{2s_n}f_2\cdot\ldots\cdot T^{ks_n}f_k)
$$
or, writing $S=\{s_n\colon n\geq 1\}$, the limit of
\begin{equation}\label{E:averages}
\frac{1}{|S\cap[0,N)|}\sum_{s\in S, s<N}( T^{s}f_1\cdot
T^{2s}f_2\cdot\ldots\cdot T^{ks}f_k)\
\end{equation}
in $L^2(\mu)$ as $N\to \infty$.

In~\cite{HK}, it is shown that the limit exists for
$s_n = n$ (another proof is given in~\cite{Z}).
If $q(n)$ is an
integer valued polynomial, the existence of the limit
for  $\{q(n)\colon n\in \mathbb{N}\}$ is shown in~\cite{HK2} (see~\cite{L} for a more
general result). It was a folklore theorem that for $k=1$, the limit exists
for the primes; it follows easily from results on
exponential sums due to Vinogradov and was written down explicitly
by Wierdl~\cite{W}. (Bourgain~\cite{Bo} and
Wierdl~\cite{W2} actually proved several stronger results on
pointwise convergence along primes.)

Let $\pi(N)$ denote the number of primes less
than or equal to $N$. Using Theorem~\ref{C:GT}
we show:
\begin{theorem}\label{T:convergence1}
Let $(X,\mathcal{X},\mu,T)$ be a  measure preserving system and
let $f_1$, $f_2\colon X\to \mathbb{C}$ be bounded measurable
functions. Then the limit
\begin{equation}
    \label{eq:two}
\lim_{N\to\infty}\frac{1}{\pi(N)} \sum_{p\in \mathbb{P}, p<N} (T^p
f_1 \cdot T^{2p}f_{2})
\end{equation}
exists in $L^2(\mu)$.
\end{theorem}

Moreover, we show that for $k=2$ a
certain factor, known as the Kronecker factor, controls the
limiting behavior  in $L^2(\mu)$ of the corresponding multiple
ergodic averages. Furthermore, in a totally ergodic system
(meaning $T$ and all its powers are ergodic), we show that the
average over the primes has the same limit as the average over the
 set of positive integers.

To prove multiple recurrence for the shifted primes and
convergence results for the primes, we compare the multiple
ergodic average related to the primes with the standard multiple
ergodic average over the full set of integers. Using estimates in
Lemma~\ref{L:lambda} and the uniformity  estimates of Green and
Tao in Theorem~\ref{C:GT}, we then show that the difference of
some modification of the two averages tends to zero. Combining
this with an ergodic version of Roth's theorem, we prove double
recurrence for the shifted primes in Section~\ref{sec:recur}. In a
similar manner, we approach the convergence questions for the
double average over the primes.  In Section~\ref{S:convergence} we
prove double convergence and in Section~\ref{sec:Kchar} we find a
characteristic factor for this average.

\section{Notation and Gowers norms}
\label{sec:nt}

By $o_{N\to\infty; a_1,\ldots,a_k}(1)$ we mean a quantity that
depends on $N,a_1,\ldots,a_k$ and for fixed $a_1,\ldots,a_k$
converges to zero as $N\to\infty$; the convergence is assumed to
be uniform with respect to all variables that are not included in
the indices.

Let
$$
\Lambda'(n)={\bf 1}_{\mathbb{P}}(n)\cdot \log{n} \ ,$$ where ${\bf
1}_{\mathbb{P}}$ denotes the indicator function of the primes. We
have
\begin{align}
\label{eq:vm} \frac 1N \sum_{0\leq
n<N}\bigl(\Lambda(n)-\Lambda'(n)\bigr) \leq & \frac 1N \sum_{0\leq
p<\sqrt N}\frac{\log N}{\log p}\cdot\log
p \\
\nonumber = & \frac{\log N} N\,\sum_{0\leq p<\sqrt N} 1\sim\frac
2{\sqrt N} \ .
\end{align}

When $f$ is a function defined on a finite set $A$, let $\mathbb{E}\left(f(n)\,\vert\, n\in A\right) =
\mathbb{E}_{n\in A} f(n)$ denote the average value of $f$ on $A$:
$$
\mathbb{E}_{n\in A}f(n) = \frac{1}{|A|}\sum_{n\in A}f(n) \ ,
$$
where by $|A|$ we mean the number of elements in $A$.  We also use
a higher dimensional version of the expectation. For example, by
$\mathbb{E}_{m,n\in A} f(n,m)$, we mean iteration of the one
variable expectation:
$$
\mathbb{E}_{m\in A}\bigl(\E_{n\in A}f(n,m)\bigr) \ .
$$

Some of our addition and averages are taken in $\Z/N\Z$ and some
in $\Z$.  In order to distinguish between the two senses, we use
the expectation notation $\E_{n\in\Z/N\Z}$ for an average in
$\Z/N\Z$ and $\frac{1}{N}\sum_{0\leq n < N}$ for an average  in
$\Z$.

If $f\colon \Z/N\Z\to \mathbb{C}$, we inductively define:
\begin{equation*}
\norm{f}_{U_1(\Z/N\Z)}=\big|\E(f(n)\,\vert\,n\in\Z/N\Z)\big|
\end{equation*}
and
\begin{equation*}
\norm{f}_{U_{d+1}(\Z/N\Z)}=\Bigl(\E(\norm{f_h\cdot
\overline{f}}_{U_d(\Z/N\Z)}^{2^d}\,\vert\,h\in\Z/N\Z)
\Bigr)^{1/2^{d+1}} \ ,
\end{equation*}
where $f_h(n) = f(n+h)$. Gowers~\cite{gowers} showed that for
$d\geq 2$ this defines a norm on $\Z/N\Z$.

\section{Some lemmas}
For studying an average over the primes, it is convenient to
replace this average with a certain weighted average over the
integers. The next lemma enables us to do this.
\begin{lemma}\label{L:replace}
If $|a_n|\leq 1$ for $n\in \mathbb{N}$, then
$$
\Big|\frac{1}{\pi(N)}\sum_{p\in\mathbb{P},p<N}a_p-
\frac{1}{N}\sum_{0\leq n<N}(\Lambda(n)\cdot
a_n)\Big|=o_{N\to\infty}(1)
$$
where $p(N)$ denotes the greatest prime number
less than or equal to $N$.
\end{lemma}
\begin{proof}
As noted in equation~\eqref{eq:vm}, we can replace $\Lambda(n)$ by
$\Lambda'(n)$ making only a small error. Then we have
\begin{align*}
\Big|\frac{1}{\pi(N)}\sum_{p\in\mathbb{P},p<N}  a_p-
\frac{1}{N}\sum_{0\leq n<N}(\Lambda'(n)\cdot a_n)\Big| & \leq
\frac{1}{N}  \sum_{p\in\mathbb{P}, p<N} |(\log{N}-\log{p})\cdot
a_p| \\ \leq \frac{1}{N}& \sum_{p\in\mathbb{P}, p<N}
\log{N}-\frac{1}{N}\sum_{0\leq n<N} \Lambda'(n) \ .
\end{align*}
Using  equation~\eqref{eq:vm} again, we have that up to a small
error term this last difference is equal to
$$
 \frac{1}{N}\sum_{p\in\mathbb{P}, p<N}
\log{N}-\frac{1}{N}\sum_{0\leq n<N} \Lambda(n).
$$
By the prime number theorem  and the well known fact that
$\Lambda$ has mean one, this difference goes to zero as
$N\to\infty$.
\end{proof}
The next Lemma is only used in the proof of Lemma~\ref{L:lambda}.
\begin{lemma}\label{L:seminorm}
For  $k\in\mathbb{N}$, let   $\theta,
\phi_0,\ldots,\phi_{k-1}\colon \Z/N\Z\to \mathbb{C}$ be functions
with $|\phi_i|\leq 1$ for $1\leq i\leq k-1$ and let $\phi_0$ be
arbitrary. Then
\begin{align*}
\Big|\E_{ m,n\in\Z/N\Z}\big(\theta(n)\cdot \phi_0(m)\cdot
\phi_1(m+n)\cdot\ldots\cdot& \phi_{k-1}(m+(k-1)n)\big)\Big|  \\
\leq &\norm{\theta}_{U_{k}(\Z/N\Z)}\cdot
\norm{\phi_0}_{L^2(\Z/N\Z)} ,
\end{align*}
where $\norm{\phi_0}_{L^2(\Z/N\Z)}=(\E_{n\in\Z/N\Z}
|\phi_0(n)|^2)^{1/2}$.
\end{lemma}

\begin{proof}
We make use of the identity
$$
\Big|\E_{n\in\Z/ N\Z}a(n)\Big|^2= \E_{n,h\in\Z/N\Z}(a(n+h)\cdot
\overline{a(n)}) \ ,
$$
which holds for $a\colon\mathbb{Z}/N\mathbb{Z}\to\mathbb{C}$. We
use induction in $k$. For $k=1$,
\begin{align*}
\Big|\E_{m,n\in\Z/N\Z}  (\theta(n)\cdot \phi_0(m))\Big|
&=\Big|\E_{n\in\Z/N\Z}\,\theta(n)\Big|\cdot \Big| \E_{m\in\Z/N\Z}\,
\phi_0(m)\Big| \\&\leq \Big|\E_{n\in\Z/N\Z}\,\theta(n)\Big|\cdot
\norm{\phi_0}_{L^2(\Z/N\Z)}\ .
\end{align*}
Suppose the statement holds for $k=l$. We show that it also holds
for $k=l+1$. Applying Cauchy-Schwarz and the previous identity we
have that
\begin{multline*}
 \Big|\E_{m,n\in\Z/N\Z}\big(\theta(n)\cdot \phi_0(m)\cdot
\phi_1(m+n)\cdot\ldots\cdot \phi_l(m+ln)\big)\Big|^2 \\
 \leq \E_{m\in\Z/N\Z}\Big|\E_{n\in\Z/N\Z}\big(\theta(n) \cdot
\phi_1(m+n)\cdot\ldots\cdot \phi_l(m+ln)\big)\Big|^2\cdot
\norm{\phi_0}_{L^2(\Z/N\Z)} \\
= \E_{m,n,h\in\Z/N\Z}\big(\theta(n+h)\cdot
\overline{\theta(n)}\cdot \phi_{0,h}'(m)\cdot
\phi_{1,h}'(m+n)\cdot\ldots\cdot \\
 \phi_{l-1,h}'(m+(l-1)n)\big)\cdot \norm{\phi_0}_{L^2(\Z/N\Z)},
\end{multline*}
where $\phi_{i,h}'(m)=\phi_{i+1}(m+h)\cdot
\overline{\phi_{i+1}(m)}$ satisfies $|\phi_{i,h}'|\leq 1$ for
$i=0,\ldots,l-1$. The last average equals
\begin{multline*} \E_{h\in\Z/
N\Z}\Bigl(\E_{m,n\in\Z/N\Z}\bigl(\theta(n+h)\cdot
\overline{\theta(n)}\cdot \phi_{0,h}'(m)\cdot
\phi_{1,h}'(m+n)\cdot\ldots  \cdot\\
\phi_{l-1,h}'(m+(l-1)n)\bigr)\Bigr) \ .
\end{multline*}
By the
induction hypothesis and the estimate
$\norm{\phi_{0,h}'}_{L^2(\Z/N\Z)}\leq 1$, the last average is
bounded by
\begin{align*}
 \E_{h\in\Z/N\Z} \norm{\theta(n+h)\cdot
\overline{\theta(n)}}_{U_{l}(\Z/N\Z)}
 &\leq \Big(\E_{h\in\Z/N\Z} \norm{\theta(n+h)\cdot
\overline{\theta(n)}}_{U_{l}(\Z/N\Z)}^{2^{l}}\Big)^{1/2^{l}} \\
& =
\norm{\theta}_{U_{l+1}(\Z/N\Z)}^2 \ .
\end{align*}
This completes the induction.
\end{proof}

We use $[x]$ to denote the greatest integer less than or equal to
$x$.

\begin{lemma}\label{L:lambda}
Let $k\geq 2$ and $N>k$  be  integers,  and $\theta\colon
\Z/N\Z\to\mathbb{R}$ be a function that is  zero for $[N/k]\leq
n<N$. Let $(X,\mathcal{X},\mu,T)$ be a measure preserving system
and $f_1,\ldots, f_{k-1}\colon X\to \mathbb{C}$ be measurable
functions with $\norm{f_i}_\infty\leq 1$ for $i=1,\ldots, k-1$.
Then
$$
\Big\Vert\frac{1}{[N/k]}\sum_{0\leq n<[N/k]} (\theta(n) \cdot
T^nf_1\cdot\ldots \cdot T^{(k-1)n}f_{k-1} )\Big\Vert_{L^2(\mu)}
\leq C_k\cdot \norm{\theta}_{U_{k}(\Z/N\Z)}\ ,
$$
for some constant $C_k > 0$.
\end{lemma}
\begin{proof}
Let $f_0\in L^\infty(\mu)$ be arbitrary. We apply
Lemma~\ref{L:seminorm} for the functions $\phi_{x,i}(n)\colon
\Z/N\Z\to \mathbb{R}$ defined by
$$
\phi_{x,0}(n)=
\begin{cases}
f_0(T^nx) &\text{ if } \  0\leq n< [N/k]\\
0 &\text{ if }\  [N/k]\leq n<N\ ,
\end{cases}
$$
 and $\phi_{x,i}(n)=f_i(T^nx)$ for $i=1,\ldots, k-1$. Note that
for $0\leq m,n<N$ we have
$$
\theta(n)\cdot \phi_{x,0}(m)\cdot \phi_{x,1}(m+n)\cdot\ldots\cdot
\phi_{x,k}(m+(k-1)n)=0
$$
except when $0\leq m,n<[N/k]$, in which case the above expression
is equal to
$$
\theta(n) \cdot f_0(T^mx) \cdot f_1(T^{m+n}x)\cdot\ldots \cdot
f_{k-1}(T^{m+(k-1)n}x) \ .
$$
Thus by Lemma~\ref{L:seminorm},
\begin{multline*}
\Big|\frac{1}{N^2} \sum_{0\leq m,n<[N/k]} \big(\theta(n) \cdot
f_0(T^mx) \cdot f_1(T^{m+n}x)\cdot\ldots \cdot
f_{k-1}(T^{m+(k-1)n}x)\big)\Big|
\\
\leq \norm{\theta}_{U_{k}(\Z/N\Z)} \cdot
\Big(\frac{1}{N}\sum_{0\leq n< [N/k]} |f_0(T^nx)|^2\Big)^{1/2}.
\end{multline*}

Integrating over $X$  and applying Cauchy-Schwarz to the right
hand side gives
\begin{multline*}
\Big| \int f_0(x)\cdot \frac{1}{[N/k]}\sum_{0\leq n<[N/k]}
(\theta(n) \cdot f_1(T^{n}x)\cdot\ldots \cdot
f_{k-1}(T^{(k-1)n}x))\ d\mu \Big|
\\ \leq
\Big(\frac{N}{[N/k]}\Big)^{3/2}\cdot  \norm{
\theta}_{U_{k}(\Z/N\Z)} \cdot \norm{f_0}_{L^2(\mu)}.
\end{multline*}
By duality, we have the advertised
estimate with $C_k=(2k)^{3/2}$.
\end{proof}

\section{Recurrence}
\label{sec:recur}
Theorem~\ref{T:recurrence1} follows immediately from
the next result, which we
prove using Theorem~\ref{C:GT} (the same statement holds for
$\mathbb{P} +1$, with the obvious modifications).

\begin{theorem}\label{T:recurrence2}
Let $(X,\mathcal{X},\mu,T)$ be a measure preserving system and
$A\in\mathcal{X}$ with $\mu(A)>0$. Then
$$
\liminf_{N\to\infty}\,  \frac{1}{\pi(N)} \sum_{p\in
\mathbb{P},p<N}\mu(A\cap T^{-(p-1)}A\cap T^{-2(p-1)}A)>0\ .
$$
\end{theorem}
\begin{proof}
By Lemma~\ref{L:replace} it suffices to show that for prime
numbers $N$ we have
$$
\liminf_{N\to\infty}\,  \frac{1}{N}\sum_{0\leq n<N}\bigl(
\Lambda(n+1)\cdot\mu(A\cap T^{-n}A\cap T^{-2n}A)\bigr) > 0\ .
$$
For this, it suffices to show that for some $w\in\mathbb{N}$ we
have
\begin{equation}\label{E:show}
\liminf_{N\to\infty}\, \frac{1}{[N/3]}\sum_{0\leq n<[N/3]}\big(
\tilde{\Lambda}_{w,N,1}(n)\cdot\mu(A\cap T^{-Wn}A\cap
T^{-2Wn}A)\big)>0\ ,
\end{equation}
where $W=\prod_{p\in\mathbb{P}, p<w}p$. We claim that
\begin{multline}\label{E:key1}
\lim_{N\to\infty} \frac{1}{[N/3]}\sum_{0\leq n<[N/3]}\big(
(\tilde{\Lambda}_{w,N,1}(n)-1)\cdot\mu(A\cap T^{-Wn}A\cap
T^{-2Wn}A)\big)
\\
=o_{w\to\infty}(1)\ .
\end{multline}
To see this, we first apply Cauchy-Schwarz to get
\begin{gather*}
\Big|\frac{1}{[N/3]}\sum_{0\leq n<[N/3]} \big(
(\tilde{\Lambda}_{w,N,1}(n)-1)\cdot
\mu(A\cap T^{Wn}A\cap T^{2Wn}A)\big)\Big| \\
\leq \Big\Vert\frac{1}{[N/3]}\sum_{0\leq n\in<[N/3]}
\big((\tilde{\Lambda}_{w,N,1}(n)-1)\cdot T^{Wn}{\bf 1}_A\cdot
T^{2Wn}{\bf 1}_A\big)\Big\Vert_{L^2(\mu)}\ .
\end{gather*}
By Lemma~\ref{L:lambda} this last term is  bounded by
$$
\norm{(\tilde{\Lambda}_{w,N,1}(n)-1)\cdot {\bf
1}_{[0,[N/3])}}_{U_3(\Z/N\Z)}\ ,
$$
and by Theorem~\ref{C:GT}, this is
$$
o_{N\to\infty;w}(1)+o_{w\to\infty}(1)\ .
$$
Letting $N\to\infty$ gives \eqref{E:key1}.

We now proceed to show \eqref{E:show}. Let $\mu(A)=\delta$. Roth's
theorem easily implies (see Theorem 2.1 in \cite{BHMP} for
details)
that for every measure preserving system $(Y,\mathcal{Y},\nu,S)$
and $B\in\mathcal{Y}$ with $\nu(B)\geq\delta$ there exists a
constant $c(\delta)>0$ such that
\begin{equation}\label{E:szemeredi}
\liminf_{N \to\infty} \frac{1}{[N/3]}\sum_{0\leq n<[N/3]}
\nu(B\cap S^{-n}B\cap S^{-2n}B)\geq c(\delta) \ .
\end{equation}
 Using \eqref{E:key1} and  applying \eqref{E:szemeredi} for the
systems $(X$, $\mathcal{X}$, $T^{W}$, $\mu)$, we have
\begin{align*}
\liminf_{N\to\infty} & \frac{1}{[N/3]} \sum_{0\leq n<[N/3]}
\big(\tilde{\Lambda}_{w,N,1}(n)
\cdot\mu(A\cap T^{-Wn}A\cap T^{-2Wn}A) \big)\\
\geq & \liminf_{N\to\infty}\, \frac{1}{[N/3]}\sum_{0\leq n<[N/3]}
\mu(A\cap T^{-Wn}A\cap
T^{-2Wn}A) +o_{w\to\infty}(1) \\
\geq & c(\delta)+o_{w\to\infty}(1)\ .
\end{align*}
Taking $w$ sufficiently large, we have \eqref{E:show}, completing
the proof.
\end{proof}

\section{$L^2$-Convergence}\label{S:convergence}
We now prove convergence of averages along the primes:
\begin{proof}[Proof of Theorem~\ref{T:convergence1}]
We can assume that $ \norm{f_i}_\infty \leq 1$ for $i=1,2$. By
Lemma~\ref{L:replace} it suffices to prove the corresponding
results for the weighted averages
\begin{equation}
\label{eq:double}
 \lim_{N\to\infty}\frac{1}{N} \sum_{0\leq n<N} (
\Lambda(n)\cdot T^n f_1 \cdot T^{2n}f_{2}) \ .
\end{equation}

For $x\in X$ let $ a_x(n)=f_1(T^nx)\cdot f_2(T^{2n}x)$. We claim
that
\begin{align}
\label{E:main} \Big\Vert\frac{1}{[WN/3]} \sum_{0\leq n<[WN/3]} &
(\Lambda(n)\cdot a_x(n))  - \\\notag
   \frac{1}{\phi(W)} & \sum_{\substack{0\leq r<W\\(r,W)=1}}\frac{1}{[N/3]} \sum_{0\leq
n<[N/3]}
a_x(Wn+r)\Big\Vert_{L^2(\mu)} \\
\notag &\ \  \ \
 =  o_{N\to\infty;w}(1) + o_{w\to\infty}(1)\ .
\end{align}
To prove \eqref{E:main} first note that
\begin{align*}
&\frac{1}{[WN/3]} \sum_{0\leq n<[(WN)/3]} (\Lambda(n)\cdot a_x(n))  \\
= &  \frac{1}{\phi(W)}\sum_{\substack{0\leq r<W\\(r,W)=1}}
\frac{1}{[N/3]}\sum_{0\leq n<[N/3]}
\Big(\frac{\phi(W)}{W}\cdot \Lambda(Wn+r)\cdot a_x(Wn+r)\Big) \\
& \qquad\qquad\qquad \qquad\qquad \qquad\qquad\qquad
\qquad\qquad\qquad \quad
+o_{N\to\infty;w}(1)  \\
= &  \frac{1}{\phi(W)}\sum_{\substack{0\leq r<W\\(r,W)=1}}
 \frac{1}{[N/3]}\sum_{0\leq n<[N/3]} \big( \tilde{\Lambda}_{w, N,r}(n) \cdot
a_x(Wn+r)\big)+o_{N\to\infty;w}(1)\ ,
\end{align*}
where the error terms are introduced because $\Lambda$ is
supported on the prime powers rather than the primes and the use
of integer parts. Hence,
\begin{align*}
&\Big\Vert \frac{1}{[WN/3]} \sum_{0\leq n<[WN/3]} (\Lambda(n)\cdot
a_x(n))- \\
& \qquad \qquad \qquad \frac{1}{\phi(W)}\sum_{\substack{0\leq
r<W\\(r,W)=1}} \frac{1}{[N/3]}\sum_{0\leq n<[N/3]]}
(a_x(Wn+r))\Big\Vert_{L^2(\mu)}\\
 \leq & \frac{1}{\phi(W)}\sum_{\substack{0\leq r<W\\(r,W)=1}}
\Big\Vert\frac{1}{[N/3]}\sum_{0\leq n<[N/3]}
\big((\tilde{\Lambda}_{w,N,r}(n)-1)
\cdot a_x(Wn+r)\big)\Big\Vert_{L^2(\mu)}\\
&\qquad\qquad\qquad \qquad\qquad \qquad\qquad\qquad
\qquad\qquad\qquad
+o_{N\to\infty;w}(1)\\
\leq & C\cdot \frac{1}{\phi(W)}\sum_{\substack{0\leq
r<W\\(r,W)=1}} \norm{{\bf 1}_{[0,N/3)}\cdot
(\tilde{\Lambda}_{w,N,r}(n)-1)}_{U_3(\Z/N\Z)} +
o_{N\to\infty;w}(1)\ .
\end{align*}
The last inequality follows by applying Lemma~\ref{L:lambda} with
$S=T^W$ and $g_i=T^{ri}f_i$, $i=1,2$. Now \eqref{E:main} follows
using Theorem~\ref{C:GT}.

We proceed to  show that the sequence
$$
A_x(N)=\E_{0\leq n<N}(\Lambda(n)\cdot a_x(n))
$$
converges in $L^2(\mu)$  by showing that it is a Cauchy sequence.
Let $\varepsilon>0$ and
$$
 B_{x,w,r}(N)=\frac{1}{[N/3]}\sum_{0\leq n<[N/3]}a_x(Wn+r)\ .
$$
Using \eqref{E:main} and the fact that $B_{x,w,r}(N)$ converges in
$L^2(\mu)$ we get that for some $W_0$, if $M,N$ are sufficiently
large  then
\begin{gather*}
 \Big\Vert A_x((W_0N)/3)-\frac{1}{\phi(W)}
\sum_{\substack{0\leq
r<W\\(r,W)=1}}(B_{x,w_0,r}(N))\Big\Vert_{L^2(\mu)}
\leq \varepsilon/3\ , \\
\norm{B_{x,w_0,r}(N)-B_{x,w_0,r}(M)}_{L^2(\mu)}\leq \varepsilon/3\
,
\end{gather*}
for all $0\leq r<W_0$ with $(r,W_0)=1$. Using this and the
triangle inequality we have that if $M,N$ are large enough then
$$
\norm{A_x((W_0N)/3)-A_x((W_0M)/3)}_{L^2(\mu)}\leq \varepsilon\ .
$$
Since
$$
A_x((W_0N)/3+i)=A_x((W_0 N)/3)+o_{N\to\infty}(1)
$$
for $0\leq i<W_0/3$, we conclude that $A_x(N)$ is Cauchy.
\end{proof}

\section{A characteristic factor for the average~\eqref{eq:two}}
\label{sec:Kchar}
 A \emph{factor} of a measure preserving system
$(X,\mathcal{X},\mu,T)$ is defined to be a  $T$-invariant
sub-$\sigma$-algebra $\mathcal{Z}$ of $\mathcal{X}$.
 The factor $\mathcal{Z}$ is
characteristic for $L^2$-convergence of the averages
in~\eqref{eq:two} if $f_1$ and $f_2$ can be replaced by their
conditional expectations $\E(f_1\,\vert\,\mathcal{Z})$  and
$\E(f_2\,\vert\,\mathcal{Z})$ without changing the value of the
limit, taken in $L^2(\mu)$.
 The \emph{Kronecker factor} $\mathcal{K}$ is
defined to be the smallest sub-$\sigma$-algebra of $\mathcal{X}$
with respect to which the  eigenfunctions of $T$ are measurable.
\begin{theorem}
\label{th:Kchar} Let $(X,\mathcal{X},\mu,T)$ be an ergodic
measure preserving system and let $f_1,f_2\colon X\to \mathbb{C}$
be bounded measurable functions. Then the Kronecker factor
$\mathcal{K}$ is characteristic for $L^2(\mu)$ convergence  of the
average in~\eqref{eq:two}. Furthermore, if $(X,\mathcal{X},\mu,T)$
is a totally ergodic system, then the limit of the average
in~\eqref{eq:two} is equal to
$$
\lim_{N\to\infty}\frac{1}{N}\sum_{0\leq n<N} (T^n f_1 \cdot
T^{2n}f_{2} )\ .
$$
\end{theorem}
\begin{proof}
We first show that the Kronecker factor is  characteristic. It
suffices to show that if either $f_1$ or $f_2$ is orthogonal to
the Kronecker factor then the average in~\eqref{eq:two} converges
to zero in $L^2(\mu)$. By $a_x(N)$, $A_x(N)$ and $B_{x,w,r}(N)$ we
denote the sequences defined in the proof of
Theorem~\ref{T:convergence1}. As it was shown in \cite{F1}, for
every $w$ and $r$ the sequence $B_{x,w,r}(N)$ converges to zero in
$L^2(\mu)$ as $N\to\infty$. Since $\lim_{N\to\infty} A_x(N)$
exists, by \eqref{E:main} we have  that
\begin{align*} \lim_{N\to\infty}A_x(N)&=\lim_{N\to\infty}A_x(WN/3)\\&=
\lim_{N\to\infty}B_{x,w,r}(N)+o_{w\to\infty}(1)=o_{w\to\infty}(1)\
,
\end{align*}
where all the limits are taken in $L^2(\mu)$. The result follows
by letting $w\to\infty$.

Next, we evaluate the limit for totally ergodic systems.  For any
such system we have for every $w,r\in\N$ and $f_1,f_2\in
L^\infty(\mu)$ that
\begin{equation}\label{E:dilation}
\lim_{N\to\infty}\frac{1}{N}\sum_{0\leq n< N} T^{Wn+r}f_1\cdot
T^{2Wn+r}f_2=\lim_{N\to\infty}\frac{1}{N}\sum_{0\leq n< N}
T^{n}f_1 \cdot T^{2n}f_2 \ ,
\end{equation}
where both limits are taken in $L^2(\mu)$. One can see this using
the formula for the limit of these averages given in
\cite{F1}. Equation~\eqref{E:dilation} gives
 $$
\lim_{N\to\infty}B_{x,w,r}(N)=\lim_{N\to\infty}
\frac{1}{N}\sum_{0\leq n<N} a_x(N)
$$
for every $w,r\in\N$. Since $A_x(N)$ converges in $L^2(\mu)$ we
conclude from \eqref{E:main} that for every $w$
\begin{equation*}
\lim_{N\to\infty}A_x(N)=\lim_{N\to\infty} A_x((WN)/3)=
\lim_{N\to\infty}\frac{1}{N}\sum_{0\leq n<N}a_x(n)
+o_{w\to\infty}(1)\ ,
\end{equation*}
where all the limits are taken in $L^2(\mu)$. The result follows
by letting $w\to\infty$.
\end{proof}

\section{Further generalizations}
\label{sec:further} It is natural to ask about higher order
recurrence for the shifted primes and higher order convergence
along the primes. The major missing ingredient is the  higher
order statements of Theorem~\ref{C:GT} with respect to the Gowers
norms. If such estimates hold, then again using a uniform version
of the multiple recurrence theorem of Furstenberg (contained
in~\cite{BHMP}), our proof carries over. For convergence, the
proof also carries over, using the general result on convergence
of linear averages in~\cite{HK}. Furthermore, the generalization
of Theorem~\ref{th:Kchar} holds, using the description of the
characteristic factors in~\cite{HK} and the higher analog of
identity~\eqref{E:dilation} in~\cite{Fr}.

With small modifications the argument used to prove
Theorems~\ref{T:recurrence1} and~\ref{T:convergence1} also
gives analogous recurrence and convergence results for two
commuting transformations.
(The proofs of these are almost identical to the proofs given and
so we omit them.)
Namely,
if $T_1$ and $T_2$ are commuting invertible  measure preserving
transformations of a probability space $(X,\mathcal{X}, \mu)$
and $A\in\mathcal{X}$ with $\mu(A)> 0$, then there exists
$p\in\mathbb{P}$ such that
$$
\mu(A\cap T_1^{-(p-1)}A\cap T_2^{-(p-1)}A) > 0 \ .
$$
The analogous statement holds with $p+1$ instead of $p-1$.
Replacing the role of the uniform version of Furstenberg's
recurrence theorem is a uniform version of the multidimensional
Szemer\'edi Theorem of Furstenberg and Katznelson~\cite{FuK}. If
one can obtain an estimate for the $k$-th Gowers norms of
$\tilde\Lambda-1$ for $k\geq 4$, analogous to the one of
Theorem~\ref{C:GT}, then as for a single transformation, one would
obtain the higher order commuting version.

Similarly, one can show that if
$f_1,f_2\in L^\infty(\mu)$, then the averages
$$
\frac{1}{\pi(N)}\sum_{p\in \mathbb{P}, p<N} (T_1^{p}f_1\cdot
T_2^{p}f_2)
$$
converge in $L^2(\mu)$ as $N\to \infty$. Again, we compare this
average with the standard ergodic average for two commuting
transformations, whose convergence was proven in~\cite{CL}.  Since
convergence is not known for the standard average of $k\geq 3$
commuting transformations, even with higher order estimates on the
Gowers norm of $\tilde{\Lambda}-1$, our proof would not
generalize.

\end{document}